\documentclass[11pt]{article}
\usepackage{amsfonts}
\usepackage{amssymb}
\usepackage{amsmath}
\def\sqr#1#2{{\vcenter{\vbox{\hrule height.#2pt
              \hbox{\vrule width.#2pt height#1pt \kern#1pt \vrule
width.#2pt}
              \hrule height.#2pt}}}}
\def\a{\alpha}
\def\b{\beta}

\def\o{\omega}

\def\3n{\negthinspace \negthinspace \negthinspace }
\def\2n{\negthinspace \negthinspace }
\def\1n{\negthinspace }

\def\ds{\displaystyle}
\def\O{\Omega}

\def\ms{\medskip}

\def\qq{\qquad}

\def\cd{\cdot}

\def\as{\hbox{\rm a.s.{ }}}

\def\Re{{\mathop{\rm Re}\,}}
\def\Im{{\mathop{\rm Im}\,}}

\def\({\Big (}
\def\){\Big )}
\def\[{\Big[}
\def\]{\Big]}

\def\={\buildrel \triangle \over =}

\def\be{\begin{equation}}
\def\bel{\begin{equation}\label}
\def\ee{\end{equation}}
\def\bea{\begin{eqnarray}}
\def\eea{\end{eqnarray}}
\def\bt{\begin{theorem}}
\def\et{\end{theorem}}
\def\bc{\begin{corollary}}
\def\ec{\end{corollary}}
\def\bl{\begin{lemma}}
\def\el{\end{lemma}}
\def\bp{\begin{proposition}}
\def\ep{\end{proposition}}
\def\br{\begin{remark}}
\def\er{\end{remark}}
\def\ba{\begin{array}}
\def\ea{\end{array}}
\def\bd{\begin{definition}}
\def\ed{\end{definition}}

\newtheorem{lemma}{Lemma}[section]
\newtheorem{remark}{Remark}[section]

\newtheorem{theorem}{Theorem}[section]
\newtheorem{corollary}{Corollary}[section]

\newtheorem{definition}{Definition}[section]
\newtheorem{proposition}{Proposition}[section]

\makeatletter
   
   \@addtoreset{equation}{section}
\makeatother

\begin{document}

\title{Calder\'on-Type Uniqueness Theorem for Stochastic Partial
Differential Equations}

\author{Xu Liu\thanks{School of Mathematics and Statistics, Northeast Normal
University, Changchun 130024, China. {\small\it e-mail:} {\small\tt
liuxu@amss.ac.cn}. \ms}
 ~~~
  and~~~
Xu Zhang\thanks{Key Laboratory of Systems and Control, Academy of
Mathematics and Systems Science, Chinese Academy of Sciences,
Beijing 100190, China; and BCAM-Basque Center for Applied
Mathematics, Bizkaia Technology Park, Building 500, E-48160, Derio,
Basque Country, Spain. {\small\it e-mail:} {\small\tt
xuzhang@amss.ac.cn}.}}
%\date{}
\maketitle

\begin{abstract}
In this Note, we present a Calder\'on-type uniqueness theorem on the
Cauchy problem of stochastic partial differential equations. To this
aim, we introduce the concept of stochastic pseudo-differential
operators, and establish their boundedness and other fundamental
properties. The proof of our uniqueness theorem is based on a new
Carleman-type estimate.

\end{abstract}

\section{Introduction and the main result}

In his remarkable paper \cite{1}, A.-P. Calder\'on established a
fundamental result on the uniqueness of the non-characteristic
Cauchy problem for general partial differential equations. One of
the main tools introduced in \cite{1} is a preliminary version of
the symbol calculation technique, which stimulated the appearance of
the theory of pseudo-differential operators (\cite{5,7}). Later,
Calder\'on's uniqueness theorem was extended to the operators with
characteristics of high multiplicity. We refer to \cite{9} and the
references cited therein for some deep results in this topic.

In recent years, one can find more and more studies on stochastic
partial differential equations (SPDEs for short). However, as far as
we know, there is no work addressing the uniqueness on the Cauchy
problem for general SPDEs. The main purpose of this Note is to
extend the classical Calder\'on's uniqueness theorem (in \cite{1})
to the stochastic setting. In order to present the key idea in the
simplest way, we do not pursue the full technical generality in this
Note. More precisely, we focus mainly on the Cauchy problem for
SPDEs in the case of simple characteristics.

Throughout this Note, $(\Omega, \mathcal{F},
\{\mathcal{F}_t\}_{t\geq 0}, P)$ is a given complete filtered
probability space, on which a $1$-dimensional standard Brownian
motion $\{w(t)\}_{t\ge 0}$ is defined.  Fix $T>0$, $n, m\in
\mathbb{N}\backslash\{0\}$,  and a neighborhood $U$ of the origin in
$\mathbb{R}^n$. Let $H$ be a Banach space. We denote by
$L_{\mathcal{F}}^\infty(0,T;H)$ the Banach space consisting of all
$H$-valued $\{\mathcal{F}_t\}_{t\ge 0}$-adapted bounded processes,
with the canonical norm; and  by
$L_{\mathcal{F}}^2(\O;C^m([0,T];H))$ the Banach space consisting of
all $H$-valued $\{\mathcal{F}_t\}_{t\ge 0}$-adapted $m$-th order
continuous differential  processes $X(\cd)$ such that
$\mathbb{E}(|X(\cd)|_{C^m([0,T];H)}^2)<\infty$, with the canonical
norm. We consider the  Cauchy problem for the following linear SPDE
of order $m$:
 \begin{eqnarray}\label{x2}
\left\{
\begin{array}{ll}
\frac{1}{i}d D^{m-1}_tu=\displaystyle\sum^{m-1}_{k=0}
\displaystyle\sum_{|\alpha|=m-k}a_\alpha(t,\o, x
)D^\alpha_xD^k_tu dt\\[3mm]
\qq\qq\qq+ \sum_{|\b|<m}\left[b_\b(t,\o, x )D^\b_{t,x}u
dt+c_\b(t,\o, x
)D^\b_{t,x}u dw(t)\right]  \quad\quad\mbox{in }(0,T)\times\Omega\times
U,\\[3mm]
u(0)=D_tu(0)=\cdots=D_t^{m-1}u(0)=0
\qq\qq\qq\qq\qq\quad\quad\quad\quad\mbox{in }\Omega\times U.
\end{array}\right.
\end{eqnarray}
In (\ref{x2}), $D_t=\frac{1}{i}\frac{\partial}{\partial t}$,
$D_{x_k}=\frac{1}{i}\frac{\partial}{\partial {x_k}}$, $i$ is the
imaginary unit, $\alpha$ and $\beta$ denote multi-indices, and
$a_\a,b_\b,c_\b\in
L_{\mathcal{F}}^\infty(\O;C^l([0,T]\times\overline{U})$ for any
$l\in \mathbb{N}$.

Denote by $\displaystyle p_m(t,\o, x, \tau,
\xi)\=\tau^m-\sum^{m-1}_{k=0} \sum_{|\alpha|=m-k}a_\alpha(t,\o, x
)\xi^\alpha\tau^k$ the symbol of the principal operator appeared in
(\ref{x2}), and by $\{\lambda_k(t,\o, x, \xi); k=1, \cdots, m\}$ the
characteristic roots of $p_m(t,\o, x, \tau, \xi)$ for any $(t, \o,
x, \xi)\in (0,T)\times\Omega\times U\times \mathbb{R}^n$, i.e.,
$p_m(t,\o, x, \lambda_k(t,\o, x, \xi), \xi)=0$. We introduce the
following hypotheses:

\vspace{1mm}

\begin{enumerate}
\item[{\bf (H1)}] All roots $\lambda_k(t,\o, x,  \xi)$ $(k=1, \cdots, m)$ are simple for any
 $(t,\o, x, \xi)\in (0,T)\times\Omega\times U\times\mathbb{R}^n$ and $|\xi|=1$;\vspace{1mm}
\item[{\bf (H2)}] For any $(t,\o, x, \xi)\in (0,T)\times\Omega\times U\times\mathbb{R}^n$ satisfying $|\xi|=1$
and any complex root $\lambda_k(t,\o, x, \xi)$, $|\mbox{Im}
\lambda_k(t,\o, x, \xi)|\geq \varepsilon$ for some positive constant
$\varepsilon$;\vspace{1mm}
\item[{\bf (H3)}] For any $(t,\o, x, \xi)\in (0,T)\times\Omega\times U\times\mathbb{R}^n$ satisfying $|\xi|=1$,
and any two distinct roots $\lambda_j(t,\o, x, \xi)$ and
$\lambda_k(t,\o, x, \xi)$, $|\lambda_j(t,\o, x, \xi)-\lambda_k(t, x,
\o, \xi)|\geq \varepsilon$ for some positive constant
$\varepsilon$.\vspace{1mm}
\end{enumerate}

\vspace{1mm}

The main result in this Note is stated as follows.

\vspace{1mm}

\begin{theorem}\label{x1}
Suppose that hypotheses (H1)--(H3) hold and $\ds u\in
\bigcap_{k_1+k_2=m-1}L_{\mathcal{F}}^2(\O;C^{k_1}([0,T];H^{k_2}(U))$
is a strong solution of equation (\ref{x2}). Then there exist a
neighborhood $V(\subset U)$ of the origin in $\mathbb{R}^n$ and a
sufficiently small $T'>0$ such that  $u$ vanishes in $(0,
T')\times\Omega\times V$.
\end{theorem}

\vspace{1mm}

We remark that the uniqueness problem in the above simple case is by
no means easy to treat. Indeed, to do this we need to introduce the
concept of stochastic pseudo-differential operators (SPDOs for
short) and to study their main properties. It is a little surprising
that the theory of stochastic pseudo-differential operators was not
available in the previous literatures although a related yet
different theory for random pseudo-differential operators was
introduced in \cite{DS}. Clearly, the study of stochastic
pseudo-differential operators has independent interest.

We refer to \cite{10} for the details of the proofs of the results
in this Note and other results in this context.

\section{Stochastic pseudo-differential operators}
As a crucial preliminary to prove Theorem \ref{x1}, in this section
we present some relevant notions and results on SPDOs. First, we
give the definition of symbol (for SPDOs). In the sequel, $l\in
\mathbb{R}$, $p\in [1, \infty]$, and $G$ is a domain in
$\mathbb{R}^n$.

\vspace{1mm}

\begin{definition}\label{x3}
If a complex-valued function $a\equiv a(t, \omega, x, \xi)$
satisfies:  1) $a(t, \omega, \cdot, \cdot)\in C^\infty(G\times
\mathbb{R}^n)$, $\forall\ (t, \omega)\in [0,T]\times \Omega$;\quad
2) $a(\cd, \cdot, x, \xi)$ is $\{\mathcal{F}_t\}_{t\geq 0}-$adapted,
$\forall\ (x, \xi)\in G\times \mathbb{R}^n$; and 3) For any two
multi-indices $\alpha$ and $\beta$, and any compact subset
$K\subseteq G$, there exists a measurable function $M_{\alpha,
\beta, K}(\cdot, \cdot)$ such that $M_{\alpha, \beta, K}(\cdot,
\cdot)\in L^p_{\mathcal{F}}(0,T)$ and $\left|
\partial^\alpha_\xi\partial^\beta_x a(t, \omega, x, \xi)\right|
\leq M_{\alpha, \beta, K}(t, \omega)(1+|\xi|)^{l-|\alpha|}$,
$\forall\ (t, \omega, x, \xi)\in [0,T]\times\Omega\times K\times
\mathbb{R}^n$, then $a$ is called a symbol of order $(l, p)$ and
denoted by $a\in S^l_p(G\times \mathbb{R}^n)$.
\end{definition}

\vspace{1mm}

We denote $\ds S^{\infty}_p(G\times
\mathbb{R}^n)=\bigcup_{m\in\mathbb{R}}S^{m}_p(G\times
\mathbb{R}^n),$ and $\ds S^{-\infty}_p(G\times
\mathbb{R}^n)=\bigcap_{m\in\mathbb{R}}S^{m}_p(G\times
\mathbb{R}^n)$. Now, for any $a\in S^l_p(G\times \mathbb{R}^n)$, we
define a SPDO $A$ (of order $(l, p)$) by
 $$
 (Au)(t, \omega, x) =\displaystyle (2\pi)^{-n}\int_{\mathbb{R}^n}e^{i
x\cdot\xi}a(t, \omega, x, \xi)\hat{u}(t, \omega, \xi)d\xi,
 $$
where $\hat u(t, \omega, \xi)=\displaystyle\int_G e^{-i
x\cdot\xi}u(t, \omega, x)dx$. We denote $A\in \mathcal{L}^l_p$.

Similar to the deterministic case, we define the amplitude, kernel,
adjoint operator, conjugate operator and properly supported SPDO.
Furthermore, we can also establish the asymptotic expansion of
symbol, etc. Here we state only the boundedness result and
invertibility of stochastic elliptic operators. For
$s\in\mathbb{R}$, put $ H^s_{comp}(G)=\{ u\in H^s(G)\ |\ u\in
\mathcal{E}'(G) \}$ and $H^s_{loc}(G)=\{ u\in \mathcal{D}'(G)\ |\
u\psi\in H^s(G)\mbox{ for any }\psi\in C^\infty_0(G) \}$.  We have
the following boundedness result.

\vspace{1mm}

\begin{theorem}\label{x7}
If $A\in \mathcal{L}_{\infty}^{l}$, then $A:
L^{p}_{\mathcal{F}}(0,T; H^s_{comp}(G))\rightarrow
L^{p}_{\mathcal{F}}(0,T; H^{s-l}_{loc}(G))$ is continuous.
\end{theorem}
\vspace{1mm}

In the rest of this section, we give the definition of elliptic
SPDOs and the invertibility result.

\vspace{1mm}

\begin{definition}\label{x8}
If the symbol $a\in S^l_p(G\times\mathbb{R}^n)$ satisfies that for
any compact subset $K\subseteq G$, there exist two positive
constants $C_K$ and $R_K$ such that $ |a(t, \omega, x, \xi)|\geq
C_K(1+|\xi|)^m$ for any $(t, \omega)\in [0,T]\times \Omega$, $x\in
K$ and $\xi\in \mathbb{R}^n$ with $|\xi|\geq R_K$,  then $A$ is
called to be elliptic.
\end{definition}

\vspace{1mm}

\begin{theorem}\label{x9}
Suppose that $A\in \mathcal{L}^l_\infty$ is elliptic. Then there
exist two operators $B_1, B_2\in \mathcal{L}^{-l}_{\infty}$, which
are the single-sided inverse for $A$ modulo smoothing operators,
namely, $B_1A-I\in \mathcal{L}^{-\infty}_\infty$ and $AB_2-I\in
\mathcal{L}^{-\infty}_\infty.$
\end{theorem}

\section{Proof of Theorem \ref{x1}}

In this section, we give a sketch of the proof of Theorem \ref{x1}.
Similar to the deterministic setting, we need to reduce (\ref{x2})
to a SPDE of order 1. For this, we denote
$A_k=\displaystyle\sum_{|\alpha|=m-k}a_\alpha(t,\o, x)D^\alpha_x$.
Denote by $a_k(t,\o, x, \xi)$ the symbol of $A_k$. Also, we
introduce a pseudo-differential operator $\Lambda^s$, whose symbol
is $(1+|\xi|^2)^{\frac{s}{2}}$. Put $M=(\Lambda^{m-1}u,
D_t\Lambda^{m-2}u, \cdots, D_t^{m-1}u)^\top$. Then, the first
equation of (\ref{x2}) can be rewritten as
 $
\frac{1}{i}dM=AMdt +fdt+Fdw(t),
 $
where $ A=\left(
\begin{array}{lccccr}
&0 &\Lambda  &\cdots &0\\
&\vdots &\vdots  &\vdots &\vdots\\
&0 &0   &\cdots &\Lambda\\
&A_m\Lambda^{1-m} &A_{m-1}\Lambda^{2-m} &\cdots &A_1
\end{array}
\right), $ $\ds f=(0, 0, \cdots, \sum_{|\b|<m}b_\b(t,\o, x
)D^\b_{t,x}u )^\top$,  and $\ds F=(0, 0, \cdots,
\sum_{|\b|<m}c_\b(t,\o, x )D^\b_{t,x}u )^\top$.

Denote by $A_0$ the pseudo-differential operator with the following
symbol
 $$ \sigma(A_0)=\left(
\begin{array}{lccccr}
&0 &|\xi|  &\cdots &0\\
&\vdots &\vdots  &\vdots &\vdots\\
&0 &0   &\cdots &|\xi|\\
&a_m|\xi|^{1-m} &a_{m-1}|\xi|^{2-m} &\cdots &a_1
\end{array}
\right).
$$
By (H1), the operator $A_0$ can be diagonalized. By (H2), if the
element of the diagonalization matrix is $A_1(t)+i B_1(t)$, then
either $B_1(t)=0$ or $B_1(t)$ is an elliptic SPDO.

By the boundedness and invertibility of elliptic SPDOs, the desired
uniqueness result follows from the following new Carleman-type
estimate for the operator ``$\frac{1}{i}dz-A_1(t)zdt-i B_1(t)zdt$".

\vspace{1mm}

\begin{proposition}\label{x11}
Suppose that $z\in
L_{\mathcal{F}}^2(\O;C([0,T];L^2(\mathbb{R}^{n})))$ is  an
$L^2(\mathbb{R}^{n})$-valued semimartingale and $z(0)=z(T)=0$ \as
Then,
 $$
 \begin{array}{rl}
&\mathbb{E}\int^T_0e^{\mu(t-T)^2}|z(t)|_{L^2(\mathbb{R}^n)}^2dt
+\frac{1}{\mu} \mathbb{E}\int^T_0e^{\mu(t-T)^2}|\mu(t-T)z(t)-
B_1(t)z(t)|^2_{L^2(\mathbb{R}^{n})}dt\\
&\leq\frac{4}{\mu}\Re
\mathbb{E}\int^T_0\int_{\mathbb{R}^n}e^{\mu(t-T)^2}\left[\frac{1}{i}dz-A_1(t)zdt-i
B_1(t)zdt\right]\cdot\overline{[i\mu(t-T)z-i B_1(t)z]}dx\\
&\quad-\frac{2}{\mu}\Im
\mathbb{E}\int^T_0\int_{\mathbb{R}^n}e^{\mu(t-T)^2}\left[\frac{1}{i}dz-A_1(t)zdt-i
B_1(t)zdt\right]\cdot\overline{(B_1(t)-B_1^*(t))z} dx
\\
&\quad-2\mathbb{E}\int^T_0\int_{\mathbb{R}^n}(t-T)e^{\mu(t-T)^2}|dz|^2dx-\frac{2}{\mu}\Re
\mathbb{E}\int^T_0e^{\mu(t-T)^2}(dz, B_1(dz))_{L^2(\mathbb{R}^n),
L^2(\mathbb{R}^n)}
\end{array}
 $$
for sufficiently small $\mu^{-1}$ and $T$, where $B_1^*$ denotes the
conjugate operator of $B_1$.
\end{proposition}

% The Acknowledgements are an un-numbered section
%\section*{Acknowledgements}
% Acknowledgements text here
\section*{Acknowledgements}
This work is partially supported by the NSF of China under grants
10901032 and 10831007, and by the National Basic Research Program of
China (973 Program) under grant 2011CB808002.

\end{document}